\documentclass[12pt,a4paper]{amsart}
\usepackage{amsmath,amsthm,amssymb}
\usepackage{graphicx}
\usepackage{tikz-cd}
\usepackage{paralist}
\usepackage{environ}
\usepackage{xcolor}
\usepackage{centernot}
\usepackage{mathtools,fancyhdr,txfonts,pxfonts}
\usepackage[utf8]{inputenc}
\usepackage[cm]{fullpage}
\usepackage{graphics}
\usepackage{amscd}
\usepackage{mathrsfs}
\usepackage{amsfonts}
\usepackage{lscape}
\usepackage{tikz}
\usetikzlibrary{matrix,calc}
\usepackage{multirow}
\usepackage{color}
\usepackage{graphicx}
\usepackage{marvosym}
\usepackage{mdframed}
\usepackage{paralist}
\usepackage{multicol}

\usepackage{pstricks,pst-text,pst-grad,pst-node,pst-3dplot,pstricks-add,pst-poly,pst-coil} 
\usepackage{pst-fun,pst-blur} 

\usepackage{hyperref}

\newtheorem{theorem}{Theorem}
\theoremstyle{plain}

\newtheorem{definition}{Definition}
\newtheorem{example}{Example}

\begin{document}

\title[Laplace Method for Determinant of cubic-matrix, of orders 2 and 3]{Laplace Method for calculate the Determinant of cubic-matrix of order 2 and order 3}

\author[Orgest ZAKA]{Orgest ZAKA}
\address{Orgest ZAKA: Department of Mathematics-Informatics, Faculty of Economy and Agribusiness, Agricultural University of Tirana, Tirana, Albania}
\email{ozaka@ubt.edu.al, gertizaka@yahoo.com}

\author[Armend Salihu]{Armend Salihu}
\address{Armend Salihu: Department of Computer Science, Faculty of Contemporary Sciences, South East European University, Tetovo, North Macedonia}
\email{ar.salihu@gmail.com}



\begin{abstract}
In this paper,  in continuation of our work, on the determinants of cubic -matrix of order 2 and order 3, we have analyzed the possibilities of developing the concept of determinant of cubic-matrix with three indexes, studying the possibility of their calculation according  the Laplace expansion method's.
We have noted that the concept of permutation expansion which is used for square determinants, as well as the concept of Laplace expansion method used for square and rectangular determinants, also can be utilized to be used for this new concept of 3D Determinants.  In this paper we proved that the Laplace expansion method's is also valid for cubic-matrix of order 2 and order 3,  these results are given clearly and with detailed proofs, they are also accompanied by illustrative examples. We also give an algorithmic presentation for the Laplace expansion method's.
\end{abstract}


\maketitle
\tableofcontents

\section{Introduction}
Based on the determinant of 2D square matrices \cite{Salihu6, Salihu7, Salihu8, ArtinM, BretscherO, Schneide-Barker, Lang}, as well as determinant of rectangular matrices \cite{Salihu1, Salihu2, Salihu3, Salihu4, Salihu5, Amiri-etal, Radic1, Radic2, MAKAREWICZetal} we have come to the idea of developing the concept of determinant of 3D cubic matrices, also in paper \cite{SalihuZaka1} we have studied and proved some basic properties related to the determinant of cubic-matrix of order 2 and 3. In this paper, we study the properties of the determinants of the cubic-matrix of order 2 and 3, related to the Laplace expansion method, our concept is based on permutation expansion method. Encouraged by geometric intuition, in this paper we are trying to give an idea and visualize the meaning of the determinants for the cubic-matrix.  Our early research mainly lies between geometry, algebra, matrix theory, etc., (see \cite{PetersZakaDyckAM}, \cite{ZakaDilauto}, \cite{ZakaFilipi2016}, \cite{FilipiZakaJusufi}, \cite{ZakaCollineations}, \cite{ZakaVertex}, \cite{ZakaThesisPhd}, \cite{ZakaPetersIso}, \cite{ZakaPetersOrder}, \cite{ZakaMohammedSF}, \cite{ZakaMohammedEndo}).

This paper is continuation of the ideas that arise based on previous researches of 3D matrix ring with element from any whatever field $F$ see \cite{zaka3DmatrixRing}, but here we study the case when the field $F$ is the field of real numbers $\mathbb{R}$ also is continuation of our research \cite{SalihuZaka1} related to the study of the properties of determinants for cubic-matrix of order 2 and 3. In this paper we follow a different method from the calculation of determinants of 3D matrix, which is studied in \cite{zaka3DGLnnp}. In contrast to the meaning of the determinant as a multi-scalar studied in \cite{zaka3DGLnnp}, in this paper we give a new definition, for the determinant of the 3D-cubic-matrix, which is a real-number.

In the papers \cite{zaka3DmatrixRing, zaka3DGLnnp}, have been studied in detail, properties for 3D-matrix, therefore, those studied properties are also valid for 3D-cubic-Matrix.

Our point in this paper is to provide a concept of determinant of 3D matrices. Our concept is based on permutation method used in regular square matrices, also based on the Laplace method which is used for calculating 2D square determinants \cite{DPoole, HERose}. 

\section{Preliminaries}

\subsection{3D Matrix}

The following is definition of 3D matrices provided by Zaka in 2017 (see \cite{zaka3DmatrixRing, zaka3DGLnnp}):
\begin{definition}
3-dimensional $m\times n\times p$ matrix will call, a matrix which has: \emph{m-horizontal layers} (analogous to m-rows), \emph{n-vertical page} (analogue with n-
columns in the usual matrices) and \emph{p-vertical layers} (p-1 of which are hidden).

The set of these matrix’s the write how:
\begin{equation}
M_{m \times n \times p} ({\rm F})=\{a_{i,j,k}|a_{i,j,k} \in F-{\rm field \hspace{6pt}  } \forall i=\overline{1,m};\hspace{6pt} j=\overline{1,n};\hspace{6pt} k=\overline{1,p}\}
\end{equation}
\end{definition}

\begin{figure}[h]
\includegraphics[width=0.4\textwidth]{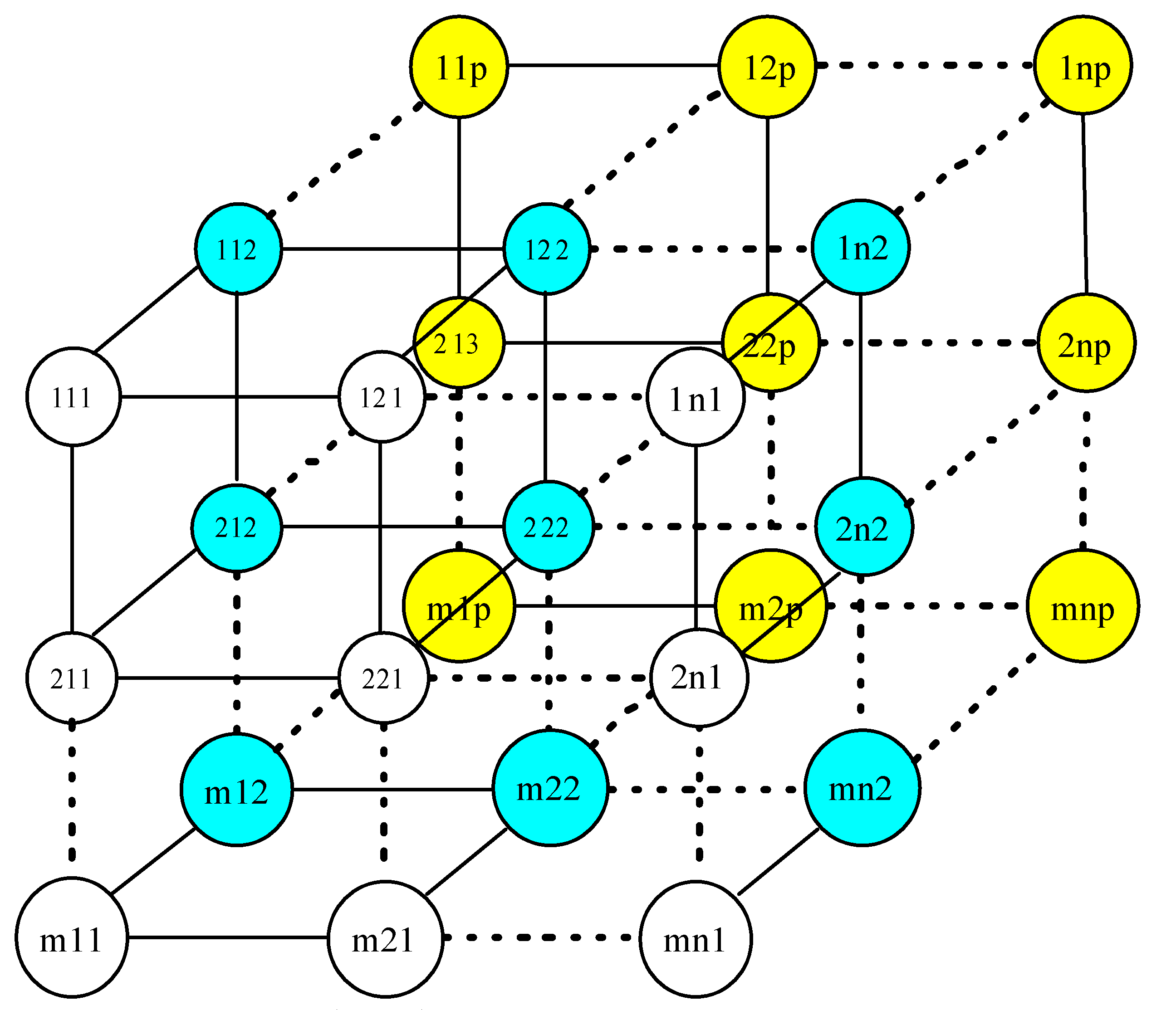}
\centering
\caption{3D-Matrix view}
\end{figure}

In the following is presented the determinant of 3D-cubic matrices, as well several properties which are adopted from 2D square determinants. 

\subsection{Cubic-Matrix of Order 2 and 3 and their Determinants}
A cubic-matrix $A_{n \times n \times n}$ for $n=2,3, \ldots$, called "cubic-matrix of order $n$".  

For $n=1$ we have that the cubic-matrix of order 1 is an element of $F$.

Let us now consider the set of cubic-matrix of order $n$, for $n=2$ or $n=3$, with elements from a field $F$ (so when cubic-matrix of order $n$, there are: $n-$vertical pages, $n-$horizontal layers and $n-$vertical layers).

From \cite{zaka3DmatrixRing, zaka3DGLnnp} we have that, the addition of 3D-matrix stands also for cubic-matrix of order 2 and 3. Also, the set of cubic-matrix of order 2 and 3 forms an commutative group (Abelian Group) related to 3Dmatrix addition.

\subsection{Determinants of Cubic-Matrix of Order 2 and 3} 

In paper \cite{SalihuZaka1}, we will define and describe the meaning of the determinants of cubic-matrix of order 2 and order 3, with elements from a field $F$.  Recall that a cubic-matrix $A_{n \times n \times n}$ for $n=2,3, \ldots$, called "cubic-matrix of order $n $".  

For $n=1$ we have that the cubic-matrix of order 1 is an element of $F$.

Let us now consider the set of cubic-matrix of order $n$, with elements from a field $F$ (so when cubic-matrix of order $n$, there are: $n-$vertical pages, $n-$horizontal layers and $n-$vertical layers),

\[
\mathcal{M}_n(F)=\{ A_{n \times n \times n}=(a_{ijk})_{n \times n \times n} |  a_{ijk} \in F, \forall i=\bar{1,n}; j=\bar{1,n}; k=\bar{1,n}  \}
\]

In this paper, we define the \emph{determinant of cubic-matrix} as a element from this field, so the map, 

\[
\begin{aligned}
\det : \mathcal{M}_n(F) & \rightarrow F \\
\forall A \in \mathcal{M}_n(F)  & \mapsto \det(A) \in F
\end{aligned}
\]

Below we give two definitions, how we will calculate the determinant of the cubic-matrix of order 2 and order 3.

\begin{definition}\label{defOrder2}
Let $A \in \mathcal{M}_2(F)$ be a $2 \times 2 \times 2$, with elements from a field $F$. 
\[
A_{2 \times 2\times 2}=
\begin{pmatrix}
\left.\begin{matrix}
a_{111} & a_{121}\\ 
a_{211} & a_{221}
\end{matrix}\right|
\begin{matrix}
a_{112} & a_{122}\\ 
a_{212} & a_{222}
\end{matrix}
\end{pmatrix}
\]

Determinant of this cubic-matrix, we called,
\[
\det[A_{2 \times 2\times 2}]=\det
\begin{pmatrix}
\left.\begin{matrix}
a_{111} & a_{121}\\ 
a_{211} & a_{221}
\end{matrix}\right|
\begin{matrix}
a_{112} & a_{122}\\ 
a_{212} & a_{222}
\end{matrix}
\end{pmatrix}=a_{111}\cdot a_{222} - a_{112}\cdot a_{221} - a_{121}\cdot a_{212} + a_{122}\cdot a_{211}
\]
\end{definition}

The follow example is case where cubic-matrix, is with elements from the number field $\mathbb{R}$.
\begin{example} \label{Example1}
Let's have the cubic-matrix, with element in number field  $\mathbb{R}$, 
\[
\det[A_{2 \times 2\times 2}]=\det
\begin{pmatrix}
\left.\begin{matrix}
4 & -3\\ 
-1 & 5
\end{matrix}\right|
\begin{matrix}
-2 & 4\\ 
-7 & 3
\end{matrix}
\end{pmatrix}
\]

then according to the definition \ref{defOrder2}, we calculate the Determinant of this cubic-matrix, and have,
\[
\det[A_{2 \times 2\times 2}]=\det
\begin{pmatrix}
\left.\begin{matrix}
4 & -3\\ 
-1 & 5
\end{matrix}\right|
\begin{matrix}
-2 & 4\\ 
-7 & 3
\end{matrix}
\end{pmatrix}
\]
\[
=4\cdot 3 - (-2)\cdot 5 - (-3)\cdot (-7) + 4\cdot (-1) = 12 - (-10) - 21 + (-4) = 12 + 10 - 21 - 4 = -3.
\]

\end{example}

We are trying to expand the meaning of the determinant of cubic-matrix, for order 3 (so when cubic-matrix, there are: 3-vertical pages, 3-horizontal layers and 3-vertical layers).

\begin{definition} \label{defOrder3}
 Let $A \in \mathcal{M}_3(F)$ be a $3 \times 3 \times 3$ cubic-matrix with element from a field $F$, 
\[
A_{3 \times 3\times 3}=
\begin{pmatrix} \left.
\begin{matrix}
a_{111} & a_{121} & a_{131}\\ 
a_{211} & a_{221} & a_{231}\\ 
a_{311} & a_{321} & a_{331}
\end{matrix}\right|
\begin{matrix}
a_{112} & a_{122} & a_{132}\\ 
a_{212} & a_{222} & a_{232}\\ 
a_{312} & a_{322} & a_{332}
\end{matrix}\left|
\begin{matrix}
a_{113} & a_{123} & a_{133}\\ 
a_{213} & a_{223} & a_{233}\\ 
a_{313} & a_{323} & a_{333}
\end{matrix}\right.
\end{pmatrix} 
. \]

Determinant of this cubic-matrix, we called,

\begin{equation}
\det[A_{3 \times 3\times 3}]=
\det\begin{pmatrix} \left.
\begin{matrix}
a_{111} & a_{121} & a_{131}\\ 
a_{211} & a_{221} & a_{231}\\ 
a_{311} & a_{321} & a_{331}
\end{matrix}\right|
\begin{matrix}
a_{112} & a_{122} & a_{132}\\ 
a_{212} & a_{222} & a_{232}\\ 
a_{312} & a_{322} & a_{332}
\end{matrix}\left|
\begin{matrix}
a_{113} & a_{123} & a_{133}\\ 
a_{213} & a_{223} & a_{233}\\ 
a_{313} & a_{323} & a_{333}
\end{matrix}\right.
\end{pmatrix} 
\end{equation}
$$=a_{111} \cdot a_{222} \cdot a_{333} - a_{111} \cdot a_{232} \cdot a_{323} - a_{111} \cdot a_{223} \cdot a_{332} + a_{111} \cdot a_{233} \cdot a_{322}$$
$$-a_{112} \cdot a_{221} \cdot a_{333} + a_{112} \cdot a_{223} \cdot a_{331} + a_{112} \cdot a_{231} \cdot a_{323} - a_{112} \cdot a_{233} \cdot a_{321}$$
$$+a_{113} \cdot a_{221} \cdot a_{332} - a_{113} \cdot a_{222} \cdot a_{331} - a_{113} \cdot a_{231} \cdot a_{322} + a_{113} \cdot a_{232} \cdot a_{321}$$
$$-a_{121} \cdot a_{212} \cdot a_{333} + a_{121} \cdot a_{213} \cdot a_{332} + a_{121} \cdot a_{232} \cdot a_{313} - a_{121} \cdot a_{233} \cdot a_{312}$$
$$+a_{122} \cdot a_{211} \cdot a_{333} - a_{122} \cdot a_{213} \cdot a_{331} - a_{122} \cdot a_{231} \cdot a_{313} + a_{122} \cdot a_{233} \cdot a_{311}$$
$$-a_{123} \cdot a_{211} \cdot a_{332} + a_{123} \cdot a_{212} \cdot a_{331} + a_{123} \cdot a_{231} \cdot a_{312} - a_{123} \cdot a_{232} \cdot a_{311}$$
$$+a_{131} \cdot a_{212} \cdot a_{323} - a_{131} \cdot a_{213} \cdot a_{322} - a_{131} \cdot a_{222} \cdot a_{313} + a_{131} \cdot a_{223} \cdot a_{312}$$
$$-a_{132} \cdot a_{211} \cdot a_{323} + a_{132} \cdot a_{213} \cdot a_{321} + a_{132} \cdot a_{221} \cdot a_{313} - a_{132} \cdot a_{223} \cdot a_{311}$$
$$+a_{133} \cdot a_{211} \cdot a_{322} - a_{133} \cdot a_{212} \cdot a_{321} - a_{133} \cdot a_{221} \cdot a_{312} + a_{133} \cdot a_{222} \cdot a_{311}$$

\end{definition}

The follow example is case where cubic-matrix, is with elements from the number field $\mathbb{R}$.
\begin{example} \label{Example2}
Let's have the cubic-matrix of order 3, with element from number field (field of real numbers) $\mathbb{R}$,
\[
\det[A_{3 \times 3\times 3}]=
\det\begin{pmatrix} \left.
\begin{matrix}
3 & 0 & -4\\ 
2 & 5 & -1\\ 
0 & 3 & -2
\end{matrix}\right|
\begin{matrix}
-2 & 4 & 0\\ 
-3 & 0 & 3\\ 
-3 & 2 & 5
\end{matrix}\left|
\begin{matrix}
5 & 1 & 0\\ 
3 & 1 & 2\\ 
0 & 4 & 3
\end{matrix}\right.
\end{pmatrix}.
\]

Then, we calculate the Determinant of this cubic-matrix following the Definition \ref{defOrder3}, and have that,

\[
\det[A_{3 \times 3\times 3}]=
\det\begin{pmatrix} \left.
\begin{matrix}
3 & 0 & -4\\ 
2 & 5 & -1\\ 
0 & 3 & -2
\end{matrix}\right|
\begin{matrix}
-2 & 4 & 0\\ 
-3 & 0 & 3\\ 
-3 & 2 & 5
\end{matrix}\left|
\begin{matrix}
5 & 1 & 0\\ 
3 & 1 & 2\\ 
0 & 4 & 3
\end{matrix}\right.
\end{pmatrix} 
\]
$$=3 \cdot 0 \cdot 3 - 3 \cdot 3 \cdot 4 - 3 \cdot 1 \cdot 5 + 3 \cdot 2 \cdot 2 - (-2) \cdot 5 \cdot 3 + (-2) \cdot 1 \cdot (-2) + (-2) \cdot (- 1) \cdot 4$$
$$-(-2) \cdot 2 \cdot 3 + 5 \cdot 5 \cdot 5 - 5 \cdot 0 \cdot (-2) - 5 \cdot (-1) \cdot 2 + 5 \cdot 3 \cdot 3 - 0 \cdot (-3) \cdot 3 + 0 \cdot 3 \cdot 5$$
$$+0 \cdot 3 \cdot 0 - 0 \cdot 2 \cdot (-3) + 4 \cdot 2 \cdot 3 - 4 \cdot 3 \cdot (-2) - 4 \cdot (-1) \cdot 0 + 4 \cdot 2 \cdot 0 - 1 \cdot 2 \cdot 5$$
$$+1 \cdot (-3) \cdot (-2) + 1 \cdot (-1) \cdot (-3) - 1 \cdot 3 \cdot 0 + (-4) \cdot (-3) \cdot 4 - (-4) \cdot 3 \cdot 2 - (-4) \cdot 0 \cdot 0 $$
$$+(-4) \cdot 1 \cdot (-3) - 0 \cdot 2 \cdot 4 + 0 \cdot 3 \cdot 3 + 0 \cdot 5 \cdot 0 - 0 \cdot 1 \cdot 0 + 0 \cdot 2 \cdot 2 - 0 \cdot (-3) \cdot 3$$
$$ - 0 \cdot 5 \cdot (-3) + 0 \cdot 0 \cdot 0$$
$$=0 -36 - 15 + 12 + 30 + 4 + 8 + 12 + 125 + 0 + 10 + 45 + 0 + 0  + 0 + 0 + 24 + 24 + 0  + 0  - 10 + 6+ 3 - 0  + 48$$
$$+ 24 + 0  + 12 - 0+ 0  + 0 - 0 + 0 + 0  + 0 + 0 = 326$$
\end{example}

\section{Minors and Co-factors of Cubic-Matrix of order 2 and 3}

In this section we will present the meaning of Minors and co-factors for cubic-matrix of order 2 and order 3.

\subsection{Minors of cubic-matrix} 

Let us start by defining minors.
\begin{definition}
Let $A_n$ be a $n\times n \times n$ cubic-matrix (with $n \geq 2$). Denote by $A_{ijk}$ the entry of cubic-matrix $A$ at the intersection of the $i-$th \textbf{horizontal layers}, $j-$th \textbf{vertical pages} and $k-$th \textbf{vertical layers}. The minor of $A_{ijk}$ is the determinant of the sub-cubic-matrix obtained from $A$ by deleting its $i-$th \textbf{horizontal layer}, $j-$\textbf{vertical page} and $k-$\textbf{vertical layer}.
\end{definition}

We now illustrate the definition with an example.
\begin{example} \label{Example3}
Let's have the cubic-matrix of order 3, with element from number field (field of real numbers) $\mathbb{R}$,
\[
A_{3 \times 3\times 3}=
\begin{pmatrix} \left.
\begin{matrix}
3 & 0 & -4\\ 
2 & 5 & -1\\ 
0 & 3 & -2
\end{matrix}\right|
\begin{matrix}
-2 & 4 & 0\\ 
-3 & 0 & 3\\ 
-3 & 2 & 5
\end{matrix}\left|
\begin{matrix}
5 & 1 & 0\\ 
3 & 1 & 2\\ 
0 & 4 & 3
\end{matrix}\right.
\end{pmatrix} 
.  \]

Take the entry  $A_{111}=3$, The sub-cubic-matrix obtained by deleting the first-horizontal layer, first-vertical page and first-vertical layer is,
\[
\begin{pmatrix}
\left.
\begin{matrix}
0 & 3\\ 
2 & 5
\end{matrix} \right |
\begin{matrix}
1 & 2\\ 
4 & 3
\end{matrix} 
\end{pmatrix}
.  \]

Thus, the minor of $A_{111}$ is

\[
M_{111}=\det\begin{pmatrix}
\left.
\begin{matrix}
0 & 3\\ 
2 & 5
\end{matrix} \right |
\begin{matrix}
1 & 2\\ 
4 & 3
\end{matrix}
\end{pmatrix}
=0 \cdot 3 - 1\cdot 5 -3 \cdot 4 + 2 \cdot 2 = -5 - 12 + 4 = -13
.  \]

Take the entry  $A_{123}=1$, The sub-cubic-matrix obtained by deleting the first-horizontal layer, 2-vertical page and 3-vertical layer is,
\[
\begin{pmatrix}
\left.
\begin{matrix}
2 &  -1\\ 
0 &  -2
\end{matrix}\right |
\begin{matrix}
-3 & 3\\ 
-3 & 5
\end{matrix}

\end{pmatrix}
.  \]
Thus, the minor of $A_{123}$ is
\[
M_{123}=\det \begin{pmatrix}
\left.
\begin{matrix}
2 &  -1\\ 
0 &  -2
\end{matrix}\right |
\begin{matrix}
-3 & 3\\ 
-3 & 5
\end{matrix}
\end{pmatrix}
=2 \cdot 5 - (-3) \cdot (-2) - (-1) \cdot (-3) + 3 \cdot 0 = 10 - 6 - 3 + 0 = 1
.  \]

\end{example}

\subsection{Co-factors of cubic-matrix of order 2 and 3} 

A co-factor is a minor whose sign may have been changed depending on the location of the respective matrix entry.
\begin{definition}\label{cofactor}

Let $A_n$ be a $n\times n \times n$ cubic-matrix (with $n \geq 2$). Denote by $M_{ijk}$ the minor of an entry $A_{ijk}$.  The co-factor of $A_{ijk}$ is
\[
C_{ijk}=(-1)^{i+j+k}\cdot M_{ijk}.
\]
\end{definition}

As an example, the pattern of sign changes $(-1)^{i+j+k}$ of a cubic-matrix of order 3 is

\[
\begin{pmatrix}
\left.
\begin{matrix}
-&+&-\\+&-&+\\-&+&-
\end{matrix}\right |
\begin{matrix}
+&-&+\\-&+&-\\+&-&+
\end{matrix}\left |
\begin{matrix}
-&+&-\\+&-&+\\-&+&-
\end{matrix} \right.
\end{pmatrix}
.  \]

\begin{example} \label{Example4}
Let's have the cubic-matrix of order 3, with element from number field (field of real numbers) $\mathbb{R}$,
\[
A_{3 \times 3\times 3}=
\begin{pmatrix} \left.
\begin{matrix}
3 & 0 & -4\\ 
2 & 5 & -1\\ 
0 & 3 & -2
\end{matrix}\right|
\begin{matrix}
-2 & 4 & 0\\ 
-3 & 0 & 3\\ 
-3 & 2 & 5
\end{matrix}\left|
\begin{matrix}
5 & 1 & 0\\ 
3 & 1 & 2\\ 
0 & 4 & 3
\end{matrix}\right.
\end{pmatrix} 
.  \]

Take the entry  $A_{111}=3$. The minor of $A_{111}$ is

\[
M_{111}=\det\begin{pmatrix}
\left.
\begin{matrix}
0 & 3\\ 
2 & 5
\end{matrix} \right |
\begin{matrix}
1 & 2\\ 
4 & 3
\end{matrix}
\end{pmatrix} = -13
\]
and its cofactor is
\[
C_{111}=(-1)^{1+1+1}\cdot M_{111}=- M_{111} = - (-13) = 13.
\]

Take the entry  $A_{123}=1$.  Thus, the minor of $A_{123}$ is
\[
M_{123}=\det \begin{pmatrix}
\left.
\begin{matrix}
2 &  -1\\ 
0 &  -2
\end{matrix}\right |
\begin{matrix}
-3 & 3\\ 
-3 & 5
\end{matrix}
\end{pmatrix} = 1
 \]
and its co-factor is
\[
C_{123}=(-1)^{1+2+3}\cdot M_{123}=M_{123} = 1.
\]

\end{example}

\section{Laplace Expansion for determinants of cubic-matrix of order 2 and 3} 
We are now ready to present the Laplace expansion.

Following the Laplace expansion method for 2D square-matrix, we are conjecturing this method for 3D cubic-matrix,

\begin{mdframed}[backgroundcolor=green!15]\label{Laplace} \textbf{Laplace Expansion}\\
If we have $A$ a cubic-matrix of order $2$ or $3$. Denote by $C_{ijk}$ the co-factor of an entry $A_{ijk}$. Then:
\begin{description}
\item[$L_1$]  For any horizontal layer $i$, the following 'horizontal layer' expansion holds:
\[
\det(A)=\sum_{jk} {A_{ijk}\cdot C_{ijk}}.
\]
\item[$L_2$] For any  'vertical page' $j$ , the following 'vertical page' expansion holds:
\[
\det(A)=\sum_{ik} {A_{ijk}\cdot C_{ijk}}.
\]
\item[$L_3$] For any  'vertical layer' $k$,  the following 'vertical page' expansion holds:
\[
\det(A)=\sum_{ij} {A_{ijk}\cdot C_{ijk}}.
\]
\end{description}
\end{mdframed}
\subsection{Laplace Expansion for determinants of cubic-matrix of order 2}
Below we prove that this method is valid for calculating the determinants of the cubic-matrix of order 2.

\begin{theorem} \label{LaplaceOrder2}
Let $A$ be a cubic-matrix of order $2$,
\[
A=
.
\]
The determinant of this cubic-matrix is invariant into expansion of three "ways" to Laplace expansion.
\end{theorem}

\proof We will prove all three expansion type, $L_1,L_2,L_3$. \\
$(L_1)$:  For any horizontal layer $i$ ($i=1,2$), the following 'horizontal layer' expansion holds:
\[
\det(A)=\sum_{jk} {A_{ijk}\cdot C_{ijk}}.
\]
Really, if we take $i=1$, and we consider the meaning of the minors and co-factors for the cubic-matrix of order 2, which we described above, we have:
\begin{equation}
\begin{aligned}
\det[A_{2 \times 2\times 2}] &=
%
\\
&=a_{111}\cdot a_{222} - a_{112}\cdot a_{221} - a_{121}\cdot a_{212} + a_{122}\cdot a_{211}.
\end{aligned}
\end{equation}
We see that this result is the same as the result of Definition \ref{defOrder2}.

Now similarly we take $i=2$, and we consider the meaning of the minors and co-factors for the cubic-matrix of order 2, which we described above, we have:

\begin{equation}
\begin{aligned}
\det[A_{2 \times 2\times 2}] &=
%
 \\
&=a_{211}\cdot a_{122} - a_{221}\cdot a_{112} - a_{212}\cdot a_{121} + a_{222}\cdot a_{111}.
\end{aligned}
\end{equation}

We see that this result is the same as the result of Definition \ref{defOrder2}.
\\
$(L_2)$: for any  'vertical page' $j$ , the following 'vertical page' expansion holds:
\[
\det(A)=\sum_{ik} {A_{ijk}\cdot C_{ijk}}.
\]
Really, if we take $j=1$, and we consider the meaning of the minors and co-factors for the cubic-matrix of order 2, which we described above, we have:

\begin{equation}
\begin{aligned}
\det[A_{2 \times 2\times 2}] &=
%
\\
& =a_{111}\cdot a_{222} + a_{211}\cdot a_{122} - a_{112}\cdot a_{221} + a_{212}\cdot a_{121}.
\end{aligned}
\end{equation}

We see that this result is the same as the result of Definition \ref{defOrder2}.
\\
Now similarly we take $j=2$, and we consider the meaning of the minors and co-factors for the cubic-matrix of order 2, which we described above, we have:

\begin{equation}
\begin{aligned}
\det[A_{2 \times 2\times 2}] &=
%
\\
&=-a_{121}\cdot a_{212} - a_{221}\cdot a_{112} + a_{122}\cdot a_{211} + a_{222}\cdot a_{111}.
\end{aligned}
\end{equation}

We see that this result is the same as the result of Definition \ref{defOrder2}.
\\
$(L_3)$: For any  'vertical layer' $k$ , the following 'vertical layer' expansion holds:
\[
\det(A)=\sum_{ik} {A_{ijk}\cdot C_{ijk}}.
\]
Really, if we take $k=1$, and we consider the meaning of the minors and co-factors for the cubic-matrix of order 2, which we described above, we have:

\begin{equation}
\begin{aligned}
\det[A_{2 \times 2\times 2}] &=
%
 \\
&=a_{111}\cdot a_{222} - a_{121}\cdot a_{212} + a_{211}\cdot a_{122} - a_{221}\cdot a_{112}.
\end{aligned}
\end{equation}

We see that this result is the same as the result of Definition \ref{defOrder2}.
\\
Now similarly we take $k=2$, and we consider the meaning of the minors and co-factors for the cubic-matrix of order 2, which we described above, we have:

\begin{equation}
\begin{aligned}
\det[A_{2 \times 2\times 2}] &=
%
 \\
&=-a_{112}\cdot a_{221} + a_{122}\cdot a_{221} - a_{212}\cdot a_{121} + a_{222}\cdot a_{111}.
\end{aligned}
\end{equation}

We see that this result is the same as the result of Definition \ref{defOrder2}.

\qed

The follow example is case where cubic-matrix of second order, is with elements from the number field $\mathbb{R}$.
\begin{example} \label{Example5}
Let's have the cubic-matrix, with element in number field  $\mathbb{R}$, 
\[
A_{2 \times 2\times 2}=
%
\]

then according to the Theorem \ref{LaplaceOrder2}, we calculate the Determinant of this cubic-matrix, and have,
\[
\det[A_{2 \times 2\times 2}]=\det
%
\\
&=4\cdot 3 - (-3)\cdot (-7) - (-2)\cdot 5 + 4\cdot (-1)=-3.
\end{aligned}
\]

We see that this result is the same as the result of Example \ref{Example1}.

For $i=2$, we have:

\[
\begin{aligned}
\det[A_{2 \times 2\times 2}] &=
%
 \\
&=-1\cdot 4 - 5\cdot (-2) - (-7)\cdot (-3) + 3\cdot 4 = -3.
\end{aligned}
\]

We see that this result is the same as the result of Example \ref{Example1}.

For $j=1$, we have:

\[
\begin{aligned}
\det[A_{2 \times 2\times 2}] &=
%
\\
& =4\cdot 3 + (-1)\cdot 4 - (-2)\cdot 5 - (-7)\cdot (-3) = -3.
\end{aligned}
\]

We see that this result is the same as the result of Example \ref{Example1}.

For $j=2$, we have:

\[
\begin{aligned}
\det[A_{2 \times 2\times 2}] &=
%
\\
& =3\cdot (-7) - 5\cdot (-2) - 4\cdot (-1) + 3\cdot 4 = -3.
\end{aligned}
\]

We see that this result is the same as the result of Example \ref{Example1}.

For $k=1$, we have:

\[
\begin{aligned}
\det[A_{2 \times 2\times 2}] &=
%
 \\
&=4\cdot 3 - (-3)\cdot (-7) + (-1)\cdot 4 - 5\cdot (-2) = -3.
\end{aligned}
\]

We see that this result is the same as the result of Example \ref{Example1}.

For $k=2$, we have:

\[
\begin{aligned}
\det[A_{2 \times 2\times 2}] &=
%
 \\
&=-(-2)\cdot 5 + 4\cdot 5 - (-7)\cdot (-3) + 3\cdot 4 = -3.
\end{aligned}
\]

We see that this result is the same as the result of Example \ref{Example1}.

\end{example}

\subsection{Laplace Expansion for determinants of cubic-matrix of order 3}
Below we prove that this method is valid for calculating the determinants of the cubic-matrix of order 3.
\begin{theorem} \label{LaplaceOrder3}
Let $A$ be a cubic-matrix of order $3$,
\[
A=%
 
\]

The determinant of this cubic-matrix is invariant into expansion of three "ways" to Laplace expansion.
\end{theorem}

\proof We will prove all three expansion type, $L_1,L_2,L_3$ also for third order. \\
$(L_1)$:  For any horizontal layer $i$ ($i=1,2$), the following 'horizontal layer' expansion holds:
\[
\det(A)=\sum_{jk} {A_{ijk}\cdot C_{ijk}}.
\]
Really, if we take $i=1$, and we consider the meaning of the minors and co-factors for the cubic-matrix of order 3, which we described above, we have:

\[
A=%
 
$$

After expanding further the above determinant based on Theorem \ref{LaplaceOrder2}, we see that this result is the same as the result of Definition \ref{defOrder3}.

If we take $i=2$, and we consider the meaning of the minors and co-factors for the cubic-matrix of order 3, which we described above, we have:

\[
A=%
 
$$

After expanding further the above determinant based on Theorem \ref{LaplaceOrder2}, we see that this result is the same as the result of Definition \ref{defOrder3}.

If we take $i=3$, and we consider the meaning of the minors and co-factors for the cubic-matrix of order 3, which we described above, we have:

\[
A=%
 
$$

After expanding further the above determinant based on Theorem \ref{LaplaceOrder2}, we see that this result is the same as the result of Definition \ref{defOrder3}.

If we take $j=1$, and we consider the meaning of the minors and co-factors for the cubic-matrix of order 3, which we described above, we have:

\[
A=%
 
$$

After expanding further the above determinant based on Theorem \ref{LaplaceOrder2}, we see that this result is the same as the result of Definition \ref{defOrder3}.

If we take $j=2$, and we consider the meaning of the minors and co-factors for the cubic-matrix of order 3, which we described above, we have:

\[
A=%
 
$$

After expanding further the above determinant based on Theorem \ref{LaplaceOrder2}, we see that this result is the same as the result of Definition \ref{defOrder3}.

If we take $j=3$, and we consider the meaning of the minors and co-factors for the cubic-matrix of order 3, which we described above, we have:

\[
A=%
 
$$

After expanding further the above determinant based on Theorem \ref{LaplaceOrder2}, we see that this result is the same as the result of Definition \ref{defOrder3}.

If we take $k=1$, and we consider the meaning of the minors and co-factors for the cubic-matrix of order 3, which we described above, we have:

\[
A=%
 
$$

After expanding further the above determinant based on Theorem \ref{LaplaceOrder2}, we see that this result is the same as the result of Definition \ref{defOrder3}.

If we take $k=2$, and we consider the meaning of the minors and co-factors for the cubic-matrix of order 3, which we described above, we have:

\[
A=%
 
$$

After expanding further the above determinant based on Theorem \ref{LaplaceOrder2}, we see that this result is the same as the result of Definition \ref{defOrder3}.

If we take $k=3$, and we consider the meaning of the minors and co-factors for the cubic-matrix of order 3, which we described above, we have:

\[
A=%
 
$$

After expanding further the above determinant based on Theorem \ref{LaplaceOrder2}, we see that this result is the same as the result of Definition \ref{defOrder3}.

\qed
The follow example is case where cubic-matrix of third order, is with elements from the number field $\mathbb{R}$.
\begin{example} \label{Example6}
Let's have the cubic-matrix, with element in number field  $\mathbb{R}$, 
\[
A_{3 \times 3\times 3}=%
 
\]

then according to the theorem Theorem \ref{LaplaceOrder3}, we calculate the Determinant of this cubic-matrix, and have,
\[
\det[A_{3 \times 3\times 3}]=\det %
  = 326.
$$

After expanding further the minors of above determinant based on Theorem \ref{LaplaceOrder2}, we see that this result is the same as the result of Example  \ref{Example2}.

For $i=2$, we have:

\[
A=%
 = 326.
$$

After expanding further the minors of above determinant based on Theorem \ref{LaplaceOrder2}, we see that this result is the same as the result of Example \ref{Example2}.

For $i=3$, we have:

\[
A=%
 = 326.
$$

After expanding further the minors of above determinant based on Theorem \ref{LaplaceOrder2}, we see that this result is the same as the result of Example \ref{Example2}.

For $j=1$, we have:

\[
A=%
 = 326.
$$

After expanding further the minors of above determinant based on Theorem \ref{LaplaceOrder2}, we see that this result is the same as the result of Example \ref{Example2}.

For $j=2$, we have:

\[
A=%
 = 326.
$$

After expanding further the minors of above determinant based on Theorem \ref{LaplaceOrder2}, we see that this result is the same as the result of Example \ref{Example2}.

For $j=3$, we have:

\[
A=%
 = 326.
$$

After expanding further the minors of above determinant based on Theorem \ref{LaplaceOrder2}, we see that this result is the same as the result of Example \ref{Example2}.

For $k=1$, we have:

\[
A=%
 = 326.
$$

After expanding further the minors of above determinant based on Theorem \ref{LaplaceOrder2}, we see that this result is the same as the result of Example \ref{Example2}.

For $k=2$, we have:

\[
A=%
  = 326.
$$

After expanding further the minors of above determinant based on Theorem \ref{LaplaceOrder2}, we see that this result is the same as the result of Example \ref{Example2}.

For $k=3$, we have:

\[
A=%
 = 326.
$$

After expanding further the minors of above determinant based on Theorem \ref{LaplaceOrder2}, we see that this result is the same as the result of Example \ref{Example2}.

\end{example}

From Theorem \ref{LaplaceOrder2} and Theorem \ref{LaplaceOrder3}, we have true the following Theorem,
\begin{theorem} \label{LaplaceExpansionThm}
The Laplace Expansion for Determinant calculation, applies to cubic-matrix of order 2 and cubic matrix of order 3
\end{theorem}

\subsection{Algorithmics implementation of  Determinants for cubic-matrix of order 2 and 3}

\quad

In paper \cite{SalihuZaka1} we have presented the pseudo-code of algorithm based on the permutation expansion method as presented in Definition 1.

In the following we have also presented the pseudo-code of algorithm based on the Laplace method as presented in \ref{LaplaceExpansionThm}.

\quad

\noindent\makebox[\linewidth]{\rule{\textwidth}{1.3pt}}

\textbf{P 1:} Laplace method for determinants of cubic matrices of order 2 and 3

\noindent\makebox[\linewidth]{\rule{\textwidth}{1.3pt}}

\textbf{Step 1:} Determine the order of determinant:

\hspace{1cm} $[m,n,o] = size(A);$

\textbf{Step 2:} Checking if 3D matrix is cubic:

\hspace{0.5cm} if $m ~\sim n; m \sim= o; n \sim= o;$

\hspace{1cm} disp('A is not square, cannot calculate the determinant')

\hspace{1cm} $d = 0$;

\hspace{1cm} return

\hspace{0.5cm} end

\textbf{Step 3:} Checking if 3D matrix is higher than the 3rd order:

\hspace{0.5cm} if $m > 3;$

\hspace{1cm} disp('A is higher than the third order, hence can not be calculated.')

\hspace{1cm} $d = 0$;

\hspace{1cm} return

\hspace{0.5cm} end

\textbf{Step 4:} Initialize $d=0;$

\textbf{Step 5:} Handling base case.

\hspace{0.5cm} if $m == 1$

\hspace{1cm} $d = A;$

\hspace{1cm} return

\hspace{0.5cm} end

\textbf{Step 6:} Select which plan we shall use to expand determinant:

\hspace{0.5cm} Horizontal Layer: x1 = 1 or 2 or 3; or

\hspace{0.5cm} Vertical Layer: x2 = 1 or 2 or 3; or

\hspace{0.5cm} Vertical page: x3 = 1 or 2 or 3;

\textbf{Step 7:} Calculate 3D determinant of order 2 and 3 based on Laplace methodology:

\hspace{0.5cm} Create loop from 1 to 2 or 3 (Depending on the order of cubic matrix):

\hspace{1cm} Create loop from 1 to 2 or 3 (Depending on the order of cubic matrix):

\hspace{1.5cm} If horizontal layer is selected:

\hspace{2cm} $d=d+(-1)^{\wedge}(1+x1+i+j)\ast A(x1,i,j)\ast det\_3DLaplace(A([1:x1-1 x1+1:m],[1:i-1 i+1:n],[1:j-1 j+1:m]));$

\hspace{1.5cm} end

\hspace{1.5cm} If vertical layer is selected:

\hspace{2cm} $d=d+(-1)^{\wedge}(1+i+x2+j)\ast A(i,x2,j)\ast det\_3DLaplace(A([1:i-1 i+1:m],[1:x2-1 x2+1:n],[1:j-1 j+1:m]));$

\hspace{1.5cm} end

\hspace{1.5cm} If vertical page is selected:

\hspace{2cm} $d=d+(-1)^{\wedge}(1+i+j+x3)\ast A(j,i,x3)\ast det\_3DLaplace(A([1:i-1 i+1:m],[1:j-1 j+1:n],[1:x3-1 x3+1:m]));$

\hspace{1.5cm} end

\hspace{1cm} end

\hspace{0.5cm} end

\textbf{Step 8:} Return the result of 3D determinant.

\noindent\makebox[\linewidth]{\rule{\textwidth}{1.3pt}}

\quad

\section{Declarations} $ $\\
\textbf{Funding:} No Funding.\\
\textbf{Authors' contributions:} The contribution of the authors is equal.\\
\textbf{Data availability statements:} This manuscript does not report data.\\
\textbf{Conflict of Interest Statement:} There is no conflict of interest with any funder.

\end{document}